\documentclass[11pt]{amsart}
\usepackage{amssymb,mathrsfs}
\title{Theta Series Associated With the Schr{\"o}dinger-Weil Representation  }

\author{Jae-Hyun Yang}
\address{Department of Mathematics, Inha University, Incheon
402-751, Korea}
\email{jhyang@inha.ac.kr }

\begin{document}

\newtheorem{theorem}{Theorem}[section]
\newtheorem{lemma}{Lemma}[section]
\newtheorem{proposition}{Proposition}[section]
\newtheorem{remark}{Remark}[section]
\newtheorem{definition}{Definition}[section]

\renewcommand{\theequation}{\thesection.\arabic{equation}}
\renewcommand{\thetheorem}{\thesection.\arabic{theorem}}
\renewcommand{\thelemma}{\thesection.\arabic{lemma}}
\newcommand{\BR}{\mathbb R}
\newcommand{\BQ}{\mathbb Q}
\newcommand{\bn}{\bf n}
\def\charf {\mbox{{\text 1}\kern-.24em {\text l}}}
\newcommand{\BC}{\mathbb C}
\newcommand{\BZ}{\mathbb Z}



\begin{abstract}
{In this paper, we define the Schr{\"o}dinger-Weil representation
for the Jacobi group and construct covariant maps for the
Schr{\"o}dinger-Weil representation. Using these covariant maps,
we construct Jacobi forms with respect to an arithmetic subgroup
of the Jacobi group. }
\end{abstract}
\maketitle

\newcommand\tr{\triangleright}
\newcommand\al{\alpha}
\newcommand\be{\beta}
\newcommand\g{\gamma}
\newcommand\gh{\Cal G^J}
\newcommand\G{\Gamma}
\newcommand\de{\delta}
\newcommand\e{\epsilon}
\newcommand\z{\zeta}
\newcommand\vth{\vartheta}
\newcommand\vp{\varphi}
\newcommand\om{\omega}
\newcommand\p{\pi}
\newcommand\la{\lambda}
\newcommand\lb{\lbrace}
\newcommand\lk{\lbrack}
\newcommand\rb{\rbrace}
\newcommand\rk{\rbrack}
\newcommand\s{\sigma}
\newcommand\w{\wedge}
\newcommand\fgj{{\frak g}^J}
\newcommand\lrt{\longrightarrow}
\newcommand\lmt{\longmapsto}
\newcommand\lmk{(\lambda,\mu,\kappa)}
\newcommand\Om{\Omega}
\newcommand\ka{\kappa}
\newcommand\ba{\backslash}
\newcommand\ph{\phi}
\newcommand\M{{\Cal M}}
\newcommand\bA{\bold A}
\newcommand\bH{\bold H}

\newcommand\Hom{\text{Hom}}
\newcommand\cP{\Cal P}
\newcommand\cH{\Cal H}

\newcommand\pa{\partial}

\newcommand\pis{\pi i \sigma}
\newcommand\sd{\,\,{\vartriangleright}\kern -1.0ex{<}\,}
\newcommand\wt{\widetilde}
\newcommand\fg{\frak g}
\newcommand\fk{\frak k}
\newcommand\fp{\frak p}
\newcommand\fs{\frak s}
\newcommand\fh{\frak h}
\newcommand\Cal{\mathcal}

\newcommand\fn{{\frak n}}
\newcommand\fa{{\frak a}}
\newcommand\fm{{\frak m}}
\newcommand\fq{{\frak q}}
\newcommand\CP{{\mathcal P}_g}
\newcommand\Hgh{{\mathbb H}_g \times {\mathbb C}^{(h,g)}}
\newcommand\BD{\mathbb D}
\newcommand\BH{\mathbb H}
\newcommand\CCF{{\mathcal F}_g}
\newcommand\CM{{\mathcal M}}
\newcommand\Ggh{\Gamma_{g,h}}
\newcommand\Chg{{\mathbb C}^{(h,g)}}
\newcommand\Yd{{{\partial}\over {\partial Y}}}
\newcommand\Vd{{{\partial}\over {\partial V}}}

\newcommand\Ys{Y^{\ast}}
\newcommand\Vs{V^{\ast}}
\newcommand\LO{L_{\Omega}}
\newcommand\fac{{\frak a}_{\mathbb C}^{\ast}}

\renewcommand\th{\theta}
\renewcommand\l{\lambda}
\renewcommand\k{\kappa}
\newcommand\tg{\widetilde\gamma}
\newcommand\wmo{{\mathscr W}_{\mathcal M,\Omega}}
\newcommand\hrnm{H_\BR^{(n,m)}}
\newcommand\rmn{\BR^{(m,n)}}

%
%
\begin{section}{{\bf Introduction}}
\setcounter{equation}{0}

For a given fixed positive integer $n$, we let
$${\mathbb H}_n=\,\big\{\,\Om\in \BC^{(n,n)}\,\big|\ \Om=\,^t\Om,\ \ \ \text{Im}\,\Om>0\,\big\}$$
be the Siegel upper half plane of degree $n$ and let
$$Sp(n,\BR)=\big\{ g\in \BR^{(2n,2n)}\ \big| \ ^t\!gJ_ng= J_n\ \big\}$$
be the symplectic group of degree $n$, where $F^{(k,l)}$ denotes
the set of all $k\times l$ matrices with entries in a commutative
ring $F$ for two positive integers $k$ and $l$, $^t\!M$ denotes
the transposed matrix of a matrix $M,\ \text{Im}\,\Om$ denotes the
imaginary part of $\Om$ and
$$J_n=\begin{pmatrix} 0&I_n \\
                   -I_n&0 \\ \end{pmatrix}.$$
We see that $Sp(n,\BR)$ acts on $\BH_n$ transitively by
\begin{equation*}g\cdot \Om=(A\Om+B)(C\Om+D)^{-1}, \end{equation*}
where $g=\begin{pmatrix} A&B\\ C&D\end{pmatrix}\in Sp(n,\BR)$ and
$\Om\in \BH_n.$

For two positive integers $n$ and $m$, we consider the Heisenberg
group
$$H_{\BR}^{(n,m)}=\{\,(\l,\mu;\k)\,|\ \l,\mu\in \BR^{(m,n)},\ \k\in \BR^{(m,m)},\ \
\k+\mu\,^t\l\ \text{symmetric}\ \}$$ endowed with the following
multiplication law
$$(\l,\mu;\k)\circ (\l',\mu';\k')=(\l+\l',\mu+\mu';\k+\k'+\l\,^t\mu'-
\mu\,^t\l').$$ We let
$$G^J=Sp(n,\BR)\ltimes H_{\BR}^{(n,m)}\quad \ ( \textrm{semi-direct product})$$
be the Jacobi group endowed with the following multiplication law
$$\Big(g,(\lambda,\mu;\kappa)\Big)\cdot\Big(g',(\lambda',\mu';\kappa')\Big) =\,
\Big(gg',(\widetilde{\lambda}+\lambda',\widetilde{\mu}+ \mu';
\kappa+\kappa'+\widetilde{\lambda}\,^t\!\mu'
-\widetilde{\mu}\,^t\!\lambda')\Big)$$ with $g,g'\in Sp(n,\BR),
(\lambda,\mu;\kappa),\,(\lambda',\mu';\kappa') \in
H_{\BR}^{(n,m)}$ and
$(\widetilde{\lambda},\widetilde{\mu})=(\lambda,\mu)g'$. We let
$\G_n=Sp(n,\BZ)$ be the Siegel modular group of degree $n$. We let
$$\G^J=\G_n\ltimes H_{\BZ}^{(n,m)}$$
be the Jacobi modular group. Then we have the {\it natural action}
of $G^J$ on the Siegel-Jacobi space $\BH_{n,m}:=\BH_n\times
\BC^{(m,n)}$ defined by
\begin{equation*}\Big(g,(\lambda,\mu;\kappa)\Big)\cdot (\Om,Z)=\Big(g\!\cdot\! \Om,(Z+\lambda \Om+\mu)
(C\Om+D)^{-1}\Big), \end{equation*}
where $g=\begin{pmatrix} A&B\\
C&D\end{pmatrix} \in Sp(n,\BR),\ (\lambda,\mu; \kappa)\in
H_{\BR}^{(n,m)}$ and $(\Om,Z)\in \BH_{n,m}.$ We refer to
\cite{YJ6}-\cite{YJ12} for more details on materials related to
the Siegel-Jacobi space.

\vskip 0.2cm The Weil representation for the symplectic group was
first introduced by A. Weil in \cite{W} to reformulate Siegel's
analytic theory of quadratic forms (cf.\,\cite{Si}) in terms of
the group theoretical theory. It is well known that the Weil
representation plays a central role in the study of the
transformation behaviors of the theta series. In this paper, we
define the Schr{\"o}dinger-Weil representation for the Jacobi
group $G^J$. The aim of this paper is to construct the covariant
maps for the Schr{\"o}dinger-Weil representation, and to construct
Jacobi forms with respect to an arithmetic subgroup of $\Gamma^J$
using these covariant maps.
\vspace{0.1in}\\
\indent This paper is organized as follows. In Section 2, we
discuss the Schr{\"o}dinger representation of the Heisenberg group
$H_\BR^{(n,m)}$ associated with a symmetric nonzero real matrix of
degree $m$. In Section 3, we review the concept of a Jacobi form
briefly. In Section 4, we define the Schr{\"o}dinger-Weil
representation $\omega_\CM$ of the Jacobi group $G^J$ associated
with a symmetric positive definite matrix $\CM$ and provide some
of the actions of $\omega_\CM$ on the representation space
$L^2\big(\BR^{(m,n)}\big)$ explicitly. In Section 5, we construct
the covariant maps for the Schr{\"o}dinger-Weil representation
$\omega_\CM$. In the final section we construct Jacobi forms with
respect to an arithmetic subgroup of $\Gamma^J$ using the
covariant maps obtained in Section 5.

\vskip 0.2cm \noindent {\bf Notations\,:} \ \ We denote by $\BZ$
and $\BC$ the ring of integers, and the field of complex numbers
respectively. $\BC^{\times}$ denotes the multiplicative group of
nonzero complex numbers. $T$ denotes the multiplicative group of
complex numbers of modulus one. The symbol ``:='' means that the
expression on the right is the definition of that on the left. For
two positive integers $k$ and $l$, $F^{(k,l)}$ denotes the set of
all $k\times l$ matrices with entries in a commutative ring $F$.
For a square matrix $A\in F^{(k,k)}$ of degree $k$, $\sigma(A)$
denotes the trace of $A$. For any $M\in F^{(k,l)},\ ^t\!M$ denotes
the transposed matrix of $M$. $I_n$ denotes the identity matrix of
degree $n$. We put $i=\sqrt{-1}.$ For $z\in\BC,$ we define
$z^{1/2}=\sqrt{z}$ so that $-\pi / 2 < \ \arg (z^{1/2})\leqq
\pi/2.$ Further we put $z^{\kappa/2}=\big(z^{1/2}\big)^\kappa$ for
every $\kappa\in\BZ.$

\end{section}

\vskip 1cm

%
%

\begin{section}{{\bf The Schr{\"o}dinger Representation of $H_\BR^{(n,m)}$}}
\setcounter{equation}{0}

\vskip 0.2cm
First of all, we observe that $H_{\mathbb{R}}^{(n,m)}$ is a 2-step
nilpotent Lie group. The inverse of an element
$(\lambda,\mu;\kappa)\in H_{\mathbb {R}}^{(n,m)}$ is given by
$$(\lambda,\mu;\kappa)^{-1}=(-\lambda,-\mu;-\kappa+\lambda\,^t\!\mu-\mu\,^t\!\lambda).$$
Now we set
\begin{equation*}
[\lambda,\mu;\kappa]=(0,\mu;\kappa)\circ
(\lambda,0;0)=(\lambda,\mu;\kappa-\mu\,^t\! \lambda).
\end{equation*}
\noindent Then $H_{\mathbb {R}}^{(n,m)}$ may be regarded as a
group equipped with the following multiplication
\begin{equation*}
[\lambda,\mu;\kappa]\diamond
[\lambda_0,\mu_0;\kappa_0]=[\lambda+\lambda_0,\mu+\mu_0;
\kappa+\kappa_0+\lambda\,^t\!\mu_0+\mu_0\,^t\!\lambda].
\end{equation*}
\noindent The inverse of $[\lambda,\mu;\kappa]\in H_{\mathbb
{R}}^{(n,m)}$ is given by
$$[\lambda,\mu;\kappa]^{-1}=[-\lambda,-\mu;\kappa+\lambda\,^t\!\mu+\mu\,^t\!\lambda].$$
We set
\begin{equation*}
 L=\left\{\,[0,\mu;\kappa]\in H_{\mathbb
{R}}^{(n,m)}\,\Big| \, \mu\in \mathbb {R}^{(m,n)},\
\kappa=\,^t\!\kappa\in \mathbb {R}^{(m,m)}\ \right\}.
 \end{equation*}
\noindent Then $L$ is a commutative normal subgroup of $H_{\mathbb
{R}}^{(n,m)}$. Let ${\widehat {L}}$ be the Pontrajagin dual of
$L$, i.e., the commutative group consisting of all unitary
characters of $L$. Then ${\widehat {L}}$ is isomorphic to the
additive group $\mathbb {R}^{(m,n)}\times \text{Symm}(m,\mathbb
{R})$ via
\begin{equation*}
\langle a,{\hat a}\rangle =e^{2 \pi i\sigma({\hat
{\mu}}\,^t\!\mu+{\hat {\kappa}}\kappa)}, \ \ \ a=[0,\mu;\kappa]\in
L,\ {\hat {a}}=({\hat {\mu}},{\hat {\kappa}})\in \widehat{L},
\end{equation*}

\noindent where $\text{Symm}(m,\mathbb {R})$ denotes the space of
all symmetric $m\times m$ real matrices.\par \noindent We put
\begin{equation*}
S=\left\{\,[\lambda,0;0]\in H_{\mathbb {R}}^{(n,m)}\,\Big|\
\lambda\in \mathbb {R}^{(m,n)}\, \right\}\cong \mathbb
{R}^{(m,n)}.
\end{equation*}
\noindent Then $S$ acts on $L$ as follows:
\begin{equation*}
\alpha_{\lambda}([0,\mu;\kappa])=[0,\mu;\kappa+\lambda\,^t\!\mu+\mu\,^t\!\lambda],
\ \ \ [\lambda,0,0]\in S.
\end{equation*}
\noindent We see that the Heisenberg group $\left( H_{\mathbb
{R}}^{(n,m)}, \diamond\right)$ is isomorphic to the semi-direct
product $S\ltimes L$ of $S$ and $L$ whose multiplication is given
by
$$(\lambda,a)\cdot
(\lambda_0,a_0)=\big(\lambda+\lambda_0,a+\alpha_{\lambda}(a_0)\big),\
\ \lambda,\lambda_0\in S,\ a,a_0\in L.$$ On the other hand, $S$
acts on ${\widehat {L}}$ by
\begin{equation*}
\alpha_{\lambda}^{*}({\hat {a}})=({\hat {\mu}}+2{\hat
{\kappa}}\lambda, {\hat {\kappa}}),\ \ [\lambda,0;0]\in S,\ \
a=({\hat {\mu}},{\hat {\kappa}})\in {\widehat {L}}.
\end{equation*}
\noindent Then, we have the relation $\langle
\alpha_{\lambda}(a),{\hat {a}} \rangle=\langle
a,\alpha_{\lambda}^{*} ({\hat {a}}) \rangle$ for all $a\in L$ and
${\hat {a}}\in {\widehat {L}}.$

\smallskip

 \indent We have three types of $S$-orbits in ${\widehat
{L}}.$

\smallskip

\noindent {\scshape Type I.} Let ${\hat{\kappa}} \in
\text{Symm}(m,\mathbb {R})$ be nondegenerate. The $S$-orbit of
${\hat {a}}({\hat {\kappa}})=(0,{\hat {\kappa}}) \in {\widehat
{L}}$ is given by

\begin{equation*}
\widehat{\mathcal{O}}_{\hat{\kappa}}= \left\{(2{\hat
{\kappa}}\lambda,{\hat {\kappa}}) \in {\widehat {L}}\ \Big|\
\lambda \in \mathbb{R}^{(m,n)} \right\} \cong \mathbb {R}^{(m,n)}.
\end{equation*}

\noindent {\scshape Type II.}\ \ Let
$({\hat{\mu}},{\hat{\kappa}})\in
\mathbb{R}^{(m,n)}\times\text{Symm}(m,\mathbb {R})$ with
degenerate  ${\hat{\kappa}}\neq 0.$ Then
\begin{equation*}
\widehat{\mathcal{O}}_{(\hat{\mu}, \hat{\kappa})} = \left\{ (\hat
\mu + 2\hat{\kappa}\lambda, \hat{\kappa}) \Big|\ \lambda \in
\mathbb {R}^{(m,n)} \right\} \subsetneqq \mathbb {R}^{(m,n)}\times
\{ \hat{\kappa} \} .
\end{equation*}

\noindent {\scshape Type III.} Let $\hat{y} \in
\mathbb{R}^{(m,n)}$. The $S$-orbit ${\widehat {\mathcal {O}
}}_{\hat {y}}$ of $\hat{a}(\hat{y}) = (\hat{y} ,0)$ is given by
\begin{equation*}
{\widehat {\mathcal {O} }}_{\hat {y}}=\left\{\,({\hat
{y}},0)\,\right\}={\hat {a}} ({\hat {y}}).
\end{equation*}
\noindent We have
$$
{\widehat{L}}= \left( \bigcup_{\begin{subarray}{c} \hat{\kappa}
\in \text{Symm} (m,\mathbb{R}) \\ {\hat{\kappa}} \,\
\text{nondegenerate} \end{subarray}}
 \widehat {\mathcal{O}}_ {\hat\kappa }
 \right)
  \bigcup
\left( \bigcup_{{\hat {y}}\in \mathbb{R}^{(m,n)}}{\widehat
{\mathcal {O} }}_{\hat {y}}\right)  \bigcup  \left(
\bigcup_{\begin{subarray}{c}({\hat{\mu}},{\hat {\kappa}}) \in
\mathbb{R}^{(m,n)} \times \text{Symm}(m,\mathbb{R}) \\
\hat{\kappa} \neq 0 \,\ \text{degenerate} \end{subarray}}
 {\widehat {\mathcal {O}}}
_({\hat \mu},{\hat \kappa }) \right)
$$
\noindent as a set. The stabilizer $S_{\hat {\kappa}}$ of $S$ at
${\hat {a}}({\hat {\kappa}})=(0,{\hat {\kappa}})$ is given by
\begin{equation*}
S_{\hat {\kappa}}=\{0\}.
 \end{equation*}
\noindent And the stabilizer $S_{\hat {y}}$ of $S$ at ${\hat
{a}}({\hat {y}})= ({\hat {y}},0)$ is given by
\begin{equation*}
S_{\hat {y}}=\left\{\,[\lambda,0;0]\,\Big|\ \lambda\in \mathbb
{R}^{(m,n)}\,\right\}=S \,\cong\,\mathbb {R}^{(m,n)}.
\end{equation*}
\indent In this section, for the present being we set
$H=H_{\mathbb {R}}^{(n,m)}$ for brevity. We see that $L$ is a
closed, commutative normal subgroup of $H$. Since
$(\lambda,\mu;\kappa)=(0,\mu; \kappa+\mu\,^t\!\lambda)\circ
(\lambda,0;0)$ for $(\lambda,\mu;\kappa)\in H,$ the homogeneous
space $X=L\backslash H$ can be identified with
$\mathbb{R}^{(m,n)}$ via
$$
Lh=L\circ (\lambda,0;0)\longmapsto \lambda,\ \ \
h=(\lambda,\mu;\kappa)\in H.
$$
We observe that $H$ acts on $X$ by
\begin{equation*}
(Lh)\cdot h_0=L\,(\lambda+\lambda_0,0;0)=\lambda+\lambda_0,
\end{equation*}
\noindent where $h=(\lambda,\mu;\kappa)\in H$ and
$h_0=(\lambda_0,\mu_0;\kappa_0)\in H.$

\medskip

\indent If $h=(\lambda,\mu;\kappa)\in H$, we have
\begin{equation*}
l_h=(0, \mu; \kappa+\mu\,^t\!\lambda),\ \ \ s_h=(\lambda,0;0)
\end{equation*}
\noindent in the Mackey decomposition of $h=l_h \circ s_h$
(cf.\,\cite{M}). Thus if $h_0=(\lambda_0,\mu_0;\kappa_0)\in H,$
then we have
\begin{equation*}
s_h\circ h_0=(\lambda,0;0)\circ
(\lambda_0,\mu_0;\kappa_0)=(\lambda+\lambda_0,\mu_0;
\kappa_0+\lambda\,^t\!\mu_0)
\end{equation*}
\noindent and so
\begin{equation}
l_{s_h\circ
h_0}=\big(0,\mu_0;\kappa_0+\mu_0\,^t\!\lambda_0+\lambda\,^t\!\mu_0+
\mu_0\,^t\!\lambda\big).
\end{equation}
 \indent For a real symmetric matrix
$c=\,^tc\in \textrm{Symm}(m,\BR) $ with $c\neq 0$, we consider the
unitary character $\chi_c$ of $L$ defined by
\begin{equation}
\chi_c\left((0,\mu;\kappa)\right)=e^{\pi i \sigma(c\kappa)}\,I,\ \
\ (0,\mu;\kappa)\in L,
\end{equation}
\noindent where $I$ denotes the identity mapping. Then the
representation ${\mathscr W}_c=\text{Ind}_L^H\,\chi_c$ of $H$
induced from $\chi_c$ is realized on the Hilbert space
$H(\chi_c)=L^2\big(X,d{\dot {h}},\mathbb {C}\big) \cong
L^2\left(\mathbb{R}^{(m,n)}, d\xi\right)$ as follows. If
$h_0=(\lambda_0,\mu_0;\kappa_0)\in H$ and $x=Lh\in X$ with
$h=(\lambda,\mu;\kappa)\in H,$ we have
\begin{equation}
\left({\mathscr W}_c (h_0)f\right)(x)=\chi_c (l_{s_h\circ
h_0})\left( f(xh_0)\right),\ \ f\in H(\chi_c).
\end{equation}
\noindent It follows from (2.1) that
\begin{equation}
\left( {\mathscr W}_c (h_0)f\right)(\lambda)=e^{\pi
i\sigma\{c(\kappa_0+\mu_0\,^t\!\lambda_0+
2\lambda\,^t\!\mu_0)\}}\,f(\lambda+\lambda_0),
\end{equation}

\noindent where $h_0=(\lambda_0,\mu_0;\kappa_0)\in H$ and
$\lambda\in\BR^{(m,n)}.$ Here we identified
$x=Lh$\,(resp.\,$xh_0=Lhh_0$) with $\lambda$\,(resp.\,
$\lambda+\lambda_0$). The induced representation ${\mathscr W}_c$
is called the {\it Schr{\" {o}}dinger\ representation} of $H$
associated with $\chi_c.$ Thus ${\mathscr W}_c $ is a monomial
representation.

\medskip

\begin{theorem}
 Let $c$ be a positive definite symmetric real matrix of
degree $m$. Then the Schr{\" {o}}dinger representation ${\mathscr
W}_c $ of $H$ is irreducible. \end{theorem} \noindent{\it Proof.}\
 The proof can be found in \cite{YJ1},\ Theorem 3. \hfill$\square$

\vskip 0.2cm\noindent {\bf Remark.} We refer to
\cite{YJ1}-\cite{YJ5} for more representations of the Heisenberg
group $H_\BR^{(n,m)}$ and their related topics.

\end{section}

\vskip 1cm

%
%

\begin{section}{{\bf Jacobi Forms}}
\setcounter{equation}{0}

\vskip 0.2cm  Let $\rho$ be a rational representation of
$GL(n,\mathbb{C})$ on a finite dimensional complex vector space
$V_{\rho}.$ Let ${\mathcal M}\in \mathbb R^{(m,m)}$ be a symmetric
half-integral semi-positive definite matrix of degree $m$. Let
$C^{\infty}(\BH_{n,m},V_{\rho})$ be the algebra of all
$C^{\infty}$ functions on $\BH_{n,m}$ with values in $V_{\rho}.$
For $f\in C^{\infty}(\BH_{n,m}, V_{\rho}),$ we define
\begin{align}
  & (f|_{\rho,{\mathcal M}}[(g,(\lambda,\mu;\kappa))])(\Om,Z) \notag \\
= \,& e^{-2\,\pi\, i\,\sigma\left( {\mathcal M}(Z+\lambda
\Om+\mu)(C\Om+D)^{-1}C\,{}^t(Z+\la\,\Om\,+\,\mu)\right) } \times
e^{2\,\pi\, i\,\sigma\left( {\mathcal M}(\lambda\,
\Om\,{}^t\!\lambda\,+\,2\,\lambda\,{}^t\!Z+\,\kappa+
\mu\,{}^t\!\lambda) \right)} \\
&\times\rho(C\Om+D)^{-1}f(g\!\cdot\!\Om,(Z+\lambda
\Om+\mu)(C\Om+D)^{-1}),\notag
\end{align}
where $g=\left(\begin{matrix} A&B\\ C&D\end{matrix}\right)\in
Sp(n,\mathbb R),\ (\lambda,\mu;\kappa)\in H_{\mathbb R}^{(n,m)}$
and $(\Om,Z)\in \BH_{n,m}.$
\vspace{0.1in}\\
\noindent {\bf Definition\ 3.1.}\quad Let $\rho$ and $\mathcal M$
be as above. Let
$$H_{\mathbb Z}^{(n,m)}= \{ (\lambda,\mu;\kappa)\in H_{\mathbb R}^{(n,m)}\, \vert
\, \lambda,\mu\in \mathbb Z^{(m,n)},\ \kappa\in \mathbb
Z^{(m,m)}\,\ \}.$$ A {\it Jacobi\ form} of index $\mathcal M$ with
respect to $\rho$ on a subgroup $\Gamma$ of $\Gamma_n$ of finite
index is a holomorphic function $f\in
C^{\infty}(\BH_{n,m},V_{\rho})$ satisfying the following
conditions (A) and (B):

\smallskip

\noindent (A) \,\ $f|_{\rho,{\mathcal M}}[\tilde{\gamma}] = f$ for
all $\tilde{\gamma}\in \Gamma \ltimes H_{\mathbb Z}^{(n,m)}$.

\smallskip

\noindent (B) \,\ For each $M\in\Gamma_n$, $f|_{\rho,\CM}[M]$ has
a Fourier expansion of the following form :
$$\big( f|_{\rho,{\mathcal M}}[M]\big)(\Om,Z) = \sum\limits_{T=\,{}^tT\ge0\atop \text {half-integral}}
\sum\limits_{R\in \mathbb Z^{(n,m)}} c(T,R)\cdot e^{{ {2\pi
i}\over {\lambda_\G}}\,\sigma(T\Om)}\cdot e^{2\pi i\sigma(RZ)}$$

\ \ with a suitable $\lambda_\G\in\BZ$ and
$c(T,R)\ne 0$ only if $\left(\begin{matrix} { 1\over {\lambda_\G}}T & \frac 12R\\
\frac 12\,^t\!R&{\mathcal M}\end{matrix}\right) \geqq 0$.

\medskip

\indent If $n\geq 2,$ the condition (B) is superfluous by Koecher
principle\,(\,cf.\, \cite{Zi} Lemma 1.6). We denote by
$J_{\rho,\mathcal M}(\Gamma)$ the vector space of all Jacobi forms
of index $\mathcal{M}$ with respect to $\rho$ on $\Gamma$.
Ziegler\,(\,cf.\,\cite{Zi} Theorem 1.8 or \cite{EZ} Theorem 1.1\,)
proves that the vector space $J_{\rho,\mathcal {M}}(\Gamma)$ is
finite dimensional. In the special case $\rho(A)=(\det(A))^k$ with
$A\in GL(n,\BC)$ and a fixed $k\in\BZ$, we write $J_{k,\CM}(\G)$
instead of $J_{\rho,\CM}(\G)$ and call $k$ the {\it weight} of the
corresponding Jacobi forms. For more results on Jacobi forms with
$n>1$ and $m>1$, we refer to \cite{YJ6}-\cite{YJ9} and \cite{Zi}.

\vskip 0.2cm
 \noindent {\bf Definition\ 3.2.}\quad A Jacobi form $f\in J_{\rho,\mathcal {M}}(\Gamma)$ is said to be
a {\it cusp}\,(\,or {\it cuspidal}\,) form if $\begin{pmatrix} {
1\over {\lambda_\G}}T & {\frac 12}R\\ {\frac 12}\,^t\!R & \mathcal
{M}\end{pmatrix} > 0$ for any $T,\,R$ with $c(T,R)\ne 0.$ A Jacobi
form $f\in J_{\rho,\mathcal{M}}(\Gamma)$ is said to be {\it
singular} if it admits a Fourier expansion such that a Fourier
coefficient $c(T,R)$ vanishes unless $\text{det}\begin{pmatrix} {
1\over {\lambda_\G}}T &{\frac 12}R\\ {\frac 12}\,^t\!R & \mathcal
{M}\end{pmatrix}=0.$

\medskip We allow a weight $k$ to be half-integral.

\medskip
\noindent {\bf Definition\ 3.3.} Let $\G\subset \G_n$ be a
subgroup of finite index. A holomorphic function $f:\BH_{n,m}\lrt
\BC$ is said to be a Jacobi form of a weight $k\in {\frac 12}\BZ$
with level $\G$ and index $\CM$ if it satisfies the following
transformation formula
\begin{equation}
f({\widetilde\g}\cdot
(\Om,Z))=\,\chi(\widetilde\g)\,J_{k,\CM}({\widetilde\g},(\Om,Z))f(\Om,Z)\quad
\textrm{for\ all}\ {\widetilde\g}\in {\widetilde\G}=\G\ltimes
H_\BZ^{(n,m)},
\end{equation}

\noindent where $\chi$ is a character of ${\widetilde\G}$ and
$J_{k,\CM}:{\widetilde\G}\times\BH_{n,m}\lrt\BC^{\times}$ is an
automorphic factor defined by

\begin{eqnarray*}
J_{k,\CM}\big( {\widetilde \g},(\Om,Z)\big)=& e^{2\,  \pi\,
i\,\sigma\big({\mathcal{M}}(Z+\lambda
\Om+\mu)(C\Om+D)^{-1}C\,{}^t(Z+\lambda \Om+\mu) \big)} \hskip 3cm
\\ & \times e^{-2\pi i\sigma\left( {\mathcal{M}}(\lambda
\Om\,{}^t\!\lambda+2\lambda\,{}^t\!Z\,+\,\kappa\,+\,
\mu\,{}^t\!\lambda)\right) } \det(C\Om+D)^k\nonumber\hskip 2cm
\end{eqnarray*}
with ${\widetilde \g}=(\g,(\lambda,\mu;\kappa))\in {\widetilde\G}$ with $\g=\begin{pmatrix} A&B\\
C&D\end{pmatrix}\in \G,\ (\lambda,\mu;\kappa)\in H_{\BZ}^{(n,m)}$
and $(\Om,Z)\in \BH_{n,m}.$

\end{section}

\vskip 1cm

%
%

\begin{section}{{\bf The Schr{\"o}dinger-Weil Representation}}
\setcounter{equation}{0}

\vskip 0.2cm Throughout this section we assume that $\CM$ is a
symmetric integral positive definite $m\times m$ matrix. We
consider the Schr{\"o}dinger representation ${\mathscr W}_\CM$ of
the Heisenberg group $\hrnm$ with the central character ${\mathscr
W}_\CM((0,0;\kappa))=\chi_\CM((0,0;\kappa))=e^{\pi
i\,\s(\CM\kappa)},\ \kappa\in\text{Symm}(m,\BR)$\,(cf.\,(2.2)). We
note that the symplectic group $Sp(n,\BR)$ acts on $\hrnm$ by
conjugation inside $G^J$. For a fixed element $g\in Sp(n,\BR)$,
the irreducible unitary representation ${\mathscr W}_\CM^g$ of
$\hrnm$ defined by
\begin{equation}
{\mathscr W}_\CM^g(h)={\mathscr W}_\CM(ghg^{-1}),\quad h\in\hrnm
\end{equation}
has the property that

\begin{equation*}
{\mathscr W}_\CM^g((0,0;\k))={\mathscr W}_\CM((0,0;\k))=e^{\pi
i\,\s(\CM \k)}\,\textrm{Id}_{H(\chi_\CM)},\quad \k\in
\text{Symm}(m,\BR).
\end{equation*}
Here $\textrm{Id}_{H(\chi_\CM)}$ denotes the identity operator on
the Hilbert space $H(\chi_\CM).$ According to Stone-von Neumann
theorem, there exists a unitary operator $R_\CM(g)$ on
$H(\chi_\CM)$ such that $R_\CM(g){\mathscr W}_\CM(h)={\mathscr
W}_\CM^g(h) R_\CM(g)$ for all $h\in\hrnm.$ We observe that
$R_\CM(g)$ is determined uniquely up to a scalar of modulus one.
From now on, for brevity, we put $G=Sp(n,\BR).$ According to
Schur's lemma, we have a map $c_\CM:G\times G\lrt T$ satisfying
the relation
\begin{equation*}
R_\CM(g_1g_2)=c_\CM(g_1,g_2)R_\CM(g_1)R_\CM(g_2)\quad \textrm{for
all }\ g_1,g_2\in G.
\end{equation*}
Therefore $R_\CM$ is a projective representation of $G$ on
$H(\chi_\CM)$ and $c_\CM$ defines the cocycle class in $H^2(G,T).$
The cocycle $c_\CM$ yields the central extension $G_\CM$ of $G$ by
$T$. The group $G_\CM$ is a set $G\times T$ equipped with the
following multiplication

\begin{equation*}
(g_1,t_1)\cdot (g_2,t_2)=\big(g_1g_2,t_1t_2\,
c_\CM(g_1,g_2)^{-1}\,\big),\quad g_1,g_2\in G,\ t_1,t_2\in T.
\end{equation*}
We see immediately that the map ${\widetilde R}_\CM:G_\CM\lrt
GL(H(\chi_\CM))$ defined by

\begin{equation}
{\widetilde R}_\CM(g,t)=t\,R_\CM(g) \quad \textrm{for all}\
(g,t)\in G_\CM
\end{equation}
is a true representation of $G_\CM.$ As in Section 1.7 in
\cite{LV}, we can define the map $s_\CM:G\lrt T$ satisfying the
relation

\begin{equation*}
c_\CM(g_1,g_2)^2=s_\CM(g_1)^{-1}s_\CM(g_2)^{-1}s_\CM(g_1g_2)\quad
\textrm{for all}\ g_1,g_2\in G.
\end{equation*}
Thus we see that

\begin{equation*}
G_{2,\CM}=\left\{\, (g,t)\in G_\CM\,|\ t^2=s_\CM(g)^{-1}\,\right\}
\end{equation*}

\noindent is the metaplectic group associated with $\CM$ that is a
two-fold covering group of $G$. The restriction $R_{2,\CM}$ of
${\widetilde R}_\CM$ to $G_{2,\CM}$ is the Weil representation of
$G$ associated with $\CM$. Now we define the projective
representation $\pi_\CM$ of the Jacobi group $G^J$ by

\begin{equation}
\pi_\CM(hg)={\mathscr W}_\CM(h)\,R_\CM(g),\quad h\in\hrnm,\ g\in
G.
\end{equation}
The projective representation $\pi_\CM$ of $G^J$ is naturally
extended to the true representation $\omega_\CM$ of the group
$G_{2,\CM}^J\!=G_{2,\CM}\ltimes \hrnm.$ The representation
$\om_\CM$ is called the $\textit{Schr{\"o}dinger-Weil
representation}$ of $G^J.$ Indeed we have
\begin{equation}
\omega_\CM(h\!\cdot\!(g,t))=t\, {\mathscr
W}_\CM(h)\,R_\CM(g),\quad h\in\hrnm,\ (g,t)\in G_{2,\CM}.
\end{equation}

\vskip0.2cm We recall that the following matrices
\begin{eqnarray*}
t_0(b)&=&\begin{pmatrix} I_n& b\\
                   0& I_n\end{pmatrix}\ \textrm{with any}\
                   b=\,{}^tb\in \BR^{(n,n)},\\
g_0(\alpha)&=&\begin{pmatrix} {}^t\alpha & 0\\
                   0& \alpha^{-1}  \end{pmatrix}\ \textrm{with
                   any}\ \alpha\in GL(n,\BR),\\
\s_{n,0}&=&\begin{pmatrix} 0& -I_n\\
                   I_n&\ 0\end{pmatrix}
\end{eqnarray*}
\noindent generate the symplectic group $G=Sp(n,\BR)$
(cf.\,\cite[p.\,326]{F},\,\cite[p.\,210]{Mum}). Therefore the
following elements $h_t(\lambda,\mu;\kappa),\
t_\CM(b),\,g_\CM(\alpha)$ and $\s_{n,\CM}$ of $G_\CM\ltimes \hrnm$
defined by
\begin{eqnarray*}
&& h_t(\la,\mu;\kappa)=\big( (I_{2n},t),(\la,\mu;\kappa)\big)\
\textrm{with}\ t\in T,\ \la,\mu\in
\BR^{(m,n)}\ \textrm{and}\ \kappa\in\BR^{(m,m)} ,\\
&&t_\CM(b)=\big( (t_0(b),1),(0,0;0)\big)\ \textrm{with any}\
                   b=\,{}^tb\in \BR^{(n,n)},\\
&& g_\CM(\alpha)=\left(
\big(g_0(\alpha),1\big),(0,0;0)\right)\
\textrm{with any}\ \alpha\in GL(n,\BR),\\
 &&\s_{n,\CM}=\left( \big(\s_{n,0},1\big),(0,0;0)\right),
\end{eqnarray*}
generate the group $G_\CM\ltimes\hrnm.$ We can show that the
representation ${\widetilde R}_\CM$ is realized on the
representation $H(\chi_\CM)=L^2\big(\rmn\big)$ as follows: for
each $f\in L^2\big(\rmn\big)$ and $x\in \rmn,$ the actions of
${\widetilde R}_\CM$ on the generators are given by

\begin{equation}
\left( {\widetilde R}_\CM
\big(h_t(\lambda,\mu;\kappa)\big)f\right)(x)=\,t\,e^{\pi
i\,\s\{\CM(\kappa+\mu\,{}^t\!\lambda+2x\,{}^t\mu)\}}\,f(x+\lambda),
\end{equation}

\begin{eqnarray}
\left( {\widetilde R}_\CM\big(t_\CM(b)\big)f\right)(x)&=&e^{\pi i\,\s(\CM\, x\,b\,{}^tx)}f(x),\\
\left( {\widetilde R}_\CM\big(g_\CM(\alpha)\big)f\right)(x)&=&\big( \det \alpha\big)^{\frac m2}\,f(x\,{}^t\alpha),\\
\left( {\widetilde R}_\CM\big(\s_{n,\CM}\big)f\right)(x)&=&\left(
{\frac 1i}\right)^{\frac{mn}2} \big( \det \CM\big)^{\frac
n2}\,\int_{\rmn}f(y)\,e^{-2\,\pi i\,\s(\CM\,y\,{}^tx)}\,dy.
\end{eqnarray}

We denote by $L^2_+\big(\rmn\big)$\,$\big(
\textrm{resp.}\,\,L^2_-\big(\rmn\big)\big)$ the subspace of
$L^2\big(\rmn\big)$ consisting of even (resp.\,odd) functions in
$L^2\big(\rmn\big)$. According to Formulas (4.6)-(4.8),
$R_{2,\CM}$ is decomposed into representations of $R_{2,\CM}^\pm$

\begin{equation*}
R_{2,\CM}=R_{2,\CM}^+\oplus R_{2,\CM}^-,
\end{equation*}
where $R_{2,\CM}^+$ and $R_{2,\CM}^-$ are the even Weil
representation and the odd Weil representation of $G$ that are
realized on $L^2_+\big(\rmn\big)$ and $L^2_-\big(\rmn\big)$
respectively. Obviously the center ${\mathscr Z}^J_{2,\CM}$ of
$G_{2,\CM}^J$ is given by
\begin{equation*}
{\mathcal Z}_{2,\CM}^J=\big\{ \big( (I_{2n},1),(0,0;\k)\big)\in
G_{2,\CM}^J\,\big\} \cong \textrm{Symm}(m,\BR).
\end{equation*}
We note that the restriction of $\omega_\CM$ to $G_{2,\CM}$
coincides with $R_{2,\CM}$ and $\omega_\CM(h)={\mathscr W}_\CM(h)$
for all $h\in \hrnm.$

\vskip 0.2cm\noindent {\bf Remark.} In the case $n=m=1,\
\omega_\CM$ is dealt in \cite{BS} and \cite{Ma}. We refer to
\cite{G} and \cite{KV} for more details about the Weil
representation $R_{2,\CM}$.

\end{section}

\vskip 1cm

\newcommand\mfm{{\mathscr F}^{(\CM)} }
\newcommand\mfoz{{\mathscr F}^{(\CM)}_{\Om,Z} }
\newcommand\wg{{\widetilde g} }
\newcommand\wgm{{\widetilde \gamma} }
\newcommand\Tm{\Theta^{(\CM)} }

%
%

\begin{section}
{{\bf Covariant Maps for the Schr{\"o}dinger-Weil representation }}
\setcounter{equation}{0}

\vskip 0.2cm As before we let $\CM$ be a symmetric positive
definite $m\times m$ real matrix. We define the mapping ${\mathscr
F}^{(\CM)}:\BH_{n,m}\lrt L^2\big(\rmn\big)$ by

\begin{equation}
\mfm (\Om,Z)(x)=\,e^{\pi i\,\s\{ \CM
(x\,\Om\,{}^tx+\,2\,x\,{}^tZ)\} },\quad (\Om,Z)\in\BH_{n,m},\
x\in\rmn.
\end{equation}

\noindent For brevity we put $\mfoz:=\mfm (\Om,Z)$ for $(\Om,Z)\in
\BH_{n,m}.$

We define the automorphic factor $J_\CM:G^J\times \BH_{n,m}\lrt
\BC^{\times}$ for $G^J$ on $\BH_{n,m}$ by

\begin{eqnarray}
   J_\CM(\widetilde g,(\Om,Z))&= e^{\pi
   i\,\sigma\left({\mathcal{M}}(Z+\lambda\,
\Om+\mu)(C\Om+D)^{-1}C\,{}^t(Z+\l\,\Om+\mu)\right)} \hskip 3cm \\
& \times e^{-\pi i\,\sigma \left( \mathcal{M}(\lambda\,
\Om\,^t\!\lambda\,+\,2\,\lambda\,{}^t\!Z\,+\,\kappa\,+\,
\mu\,{}^t\!\lambda ) \right) } \det(C\Om+D)^{\frac m2},\nonumber
\end{eqnarray}
where ${\widetilde g}=(g,(\lambda,\mu;\kappa))\in G^J$ with $g=\begin{pmatrix} A&B\\
C&D\end{pmatrix}\in Sp(n,\BR),\ (\lambda,\mu;\kappa)\in
H_{\BR}^{(n,m)}$ and $(\Om,Z)\in \BH_{n,m}.$

\begin{theorem} The map ${\mathscr F}^{(\CM)}:\BH_{n,m}\lrt
L^2\big(\rmn\big)$ defined by (5.1) is a covariant map for the
Schr{\"o}dinger-Weil representation $\omega_\CM$ of $G^J$ and the
automorphic factor $J_\CM$ for $G^J$ on $\BH_{n,m}$ defined by
Formula (5.2). In other words, $\mfm$ satisfies the following
covariance relation

\begin{equation}
\om_\CM ({\widetilde g})\mfoz=J_\CM \big( {\widetilde
g},(\Om,Z)\big)^{-1} \mfm_{{\widetilde g}\cdot (\Om,Z)}
\end{equation}

\noindent for all ${\widetilde g}\in G^J$ and $(\Om,Z)\in
\BH_{n,m}.$
\end{theorem}

\noindent {\it Proof.} For an element
$\wg=(g,(\lambda,\mu;\kappa))\in G^J$ with $g=\begin{pmatrix} A &
B
\\ C & D \end{pmatrix}\in Sp(n,\BR),$ we put
$(\Om_*,Z_*)=\wg\cdot (\Om,Z)$ for $(\Om,Z)\in\BH_{n,m}.$ Then we
have
\begin{eqnarray*}
&\Om_*=g\cdot \Om=(A\Om+B)(C\Om+D)^{-1},\\
&Z_*=(Z+\l\, \Om+\mu)(C\Om+D)^{-1}.
\end{eqnarray*}

\noindent In this section we use the notations $t_0 (b),\
g_0(\alpha)$ and $\sigma_{n,0}$ in Section 4. Since the following
elements $h(\l,\mu;\kappa),\ t(b),\ g(\alpha)$ and $\s_n$ of $G^J$
defined by

\begin{eqnarray*}
 h(\l,\mu;\kappa)&=&(I_{2n},(\l,\mu;\kappa))\quad \textrm{with}\
\l,\mu\in\rmn,\ \kappa\in\BR^{(m,m)},\\
t(b)&=&\big( t_0(b),(0,0;0)\big)\quad \textrm{with}\ b=\,{}^tb\in
\BR^{(m,m)},\\
 g(\alpha)&=&\big( g_0(\alpha),(0,0;0)\big)\quad \textrm{with}\
\alpha\in GL(n,\BR),\\
 \s_n&=&\big( \s_{n,0},(0,0;0)\big)
\end{eqnarray*}

\noindent generate the Jacobi group, it suffices to prove the
covariance relation (5.3) for the above generators.

\vskip 0.5cm\noindent {\bf Case I.} $\wg=h(\l,\mu;\kappa)$ with
$\l,\mu\in\rmn,\ \kappa\in\BR^{(m,m)}.$

\vskip 0.1cm In this case, we have
$$\Om_*=\Om,\quad Z_*=Z+\la\,\Om+\mu$$
and
$$J_\CM\big(\wg,(\Om,Z)\big)=\,e^{-\pi i\,\s\{ \CM(\l\,\Om\,{}^t\l
+2\,\l\,{}^tZ+\kappa+\mu\,{}^t\!\l)\} }.$$

According to Formula (4.5), for $x\in \rmn,$

\begin{eqnarray*}
& & \left( \om_\CM \big( h(\l,\mu;\kappa)\big) \mfoz\right) (x)\\
&=& e^{\pi i\,\sigma\{ \CM(\kappa+\mu\,{}^t\l+2\,x\,{}^t\mu)\} }
\mfoz
(x+\l)\\
&=& e^{\pi i\,\sigma\{ \CM(\kappa+\mu\,{}^t\l+2\,x\,{}^t\mu)\}
}\,e^{\pi i \,\s\{ \CM(
(x+\l)\Om\,{}^t(x+\l)+\,2\,(x+\l)\,{}^tZ)\} }.
\end{eqnarray*}

\noindent On the other hand, according to Formula (5.2), for $x\in
\rmn,$

\begin{eqnarray*}
& & J_\CM \big( h(\l,\mu;\kappa),(\Om,Z)\big)^{-1}\mfm_{\wg\cdot
(\Om,Z)}(x)\\
&=& J_\CM \big(
h(\l,\mu;\kappa),(\Om,Z)\big)^{-1}\mfm_{\Om,Z+\l\,\Om+\mu}(x)\\
&=& e^{\pi i\,\sigma\{ \CM(\l
\Om\,{}^t\!\l\,+\,2\,\l\,{}^tZ\,+\,\kappa\,+\,\mu\,{}^t\!\l)\}
}\cdot e^{\pi
i\,\sigma\{ \CM (x\, \Om\,{}^tx\,+\,2\,x\,{}^t(Z+\l\Om+\mu))\} }\\
&=& e^{\pi i\,\sigma\{ \CM(\kappa+\mu\,{}^t\l+2\,x\,{}^t\mu)\}
}\,e^{\pi i \,\s\{ \CM(
(x+\l)\Om\,{}^t(x+\l)+\,2\,(x+\l)\,{}^tZ)\} }.
\end{eqnarray*}

\noindent Therefore we prove the covariance relation (5.3) in the
case $\wg=h(\l,\mu;\kappa)$ with $\l,\mu,\kappa$ real.

\vskip 0.5cm\noindent {\bf Case II.} $\wg=t(b)$ with
$b=\,{}^tb\in\BR^{(n,n)}.$

\vskip 0.1cm In this case, we have
$$\Om_*=\Om+b,\quad Z_*=Z\quad \textrm{and}\quad J_\CM\big(\wg,(\Om,Z)\big)=1.$$

\noindent According to Formula (4.6), we obtain

$$ \left( \om_\CM \big( \wg\big) \mfoz\right)
(x)=\,e^{\pi\,i\,\sigma (\CM\,xb\,{}^tx)} \mfoz (x),\quad
x\in\rmn.$$

\noindent On the other hand, according to Formula (5.2), for $x\in
\rmn,$ we obtain

\begin{eqnarray*}
& & J_\CM \big( \wg,(\Om,Z)\big)^{-1}\mfm_{\wg\cdot
(\Om,Z)}(x)\\
&=& \mfm_{\Om+b,Z}(x)\\
&=& e^{\pi i\,\sigma\left( \CM \left(
x(\Om+b)\,{}^tx+2\,x\,{}^tZ\right)\right)}\\
&=& \,e^{\pi\,i\,\sigma (\CM\,xb\,{}^tx)} \mfoz (x).
\end{eqnarray*}

\noindent Therefore we prove the covariance relation (5.3) in the
case $\wg=t(b)$ with $b=\,{}^tb\in\BR^{(n,n)}.$

\vskip 0.52cm\noindent {\bf Case III.} $\wg=g(\alpha)$ with
$\alpha\in GL(n,\BR).$

\vskip 0.1cm In this case, we have
$$\Om_*=\,{}^t\alpha\,\Om\,\alpha,\quad Z_*=Z\alpha$$
and
$$J_\CM\big(\wg,(\Om,Z)\big)=(\det\alpha)^{-{\frac m2}}.$$

\noindent According to Formula (4.7), for $x\in \rmn,$

\begin{eqnarray*}
& & \left( \om_\CM \big(\wg \big) \mfoz\right) (x)\\
&=& (\det\alpha)^{{\frac m2}} \mfoz
(x\,{}^t\alpha)\\
&=&\,(\det\alpha)^{{\frac m2}}\cdot e^{\pi i\,\sigma\{
\CM(x\,{}^t\alpha\,\Om\,{}^t(x\,{}^t\alpha)+2\,x\,{}^t\alpha\,{}^tZ)\}
}.
\end{eqnarray*}

\noindent On the other hand, according to Formula (5.2), for $x\in
\rmn,$

\begin{eqnarray*}
& & J_\CM \big( \wg,(\Om,Z)\big)^{-1}\mfm_{\wg\cdot
(\Om,Z)}(x)\\
&=&(\det\alpha)^{{\frac
m2}}\mfm_{{}^t\alpha\,\Om\,\alpha,Z\alpha}(x)\\
&=&\,(\det\alpha)^{{\frac m2}}\cdot e^{\pi i\,\sigma\{
\CM(x\,{}^t\alpha\,\Om\,{}^t(x\,{}^t\alpha)+2\,x\,{}^t\alpha\,{}^tZ)\}}.
\end{eqnarray*}

\noindent Therefore we prove the covariance relation (5.3) in the
case $\wg=g(\alpha)$ with $\alpha\in GL(n,\BR).$

\vskip 0.52cm\noindent {\bf Case IV.} $\wg=\left( \begin{pmatrix}
0 & -I_n \\ I_n & \ 0 \end{pmatrix},(0,0;0)\right).$

\vskip 0.1cm In this case, we have
$$\Om_*=-\Om^{-1},\quad Z_*=Z\,\Om^{-1}$$
and
$$J_\CM\big(\wg,(\Om,Z)\big)=e^{\pi \,i\,\sigma (\CM Z\Om^{-1}\,{}^tZ)}\,\big( \det\Om\big)^{\frac m2}.$$

\noindent In order to prove the covariance relation (5.3), we need
the following useful lemma.

\begin{lemma} For a fixed element $\Om\in \BH_n$ and a fixed
element $Z\in\BC^{(m,n)},$ we obtain the following property
\begin{equation}
\int_{\rmn} e^{\pi\,i\,\sigma
(x\,\Om\,{}^tx+2\,x\,{}^tZ)}dx_{11}\cdots dx_{mn} = \left( \det
{\Omega\over i}\right)^{-{\frac m2}}\,
e^{-\pi\,i\,\sigma(Z\,\Om^{-1}\,{}^tZ)},
\end{equation}
where $x=(x_{ij})\in \BR^{(m,n)}.$
\end{lemma}

\noindent {\it Proof of Lemma 5.1.} By a simple computation, we
see that
$$e^{\pi i\, \sigma ( x\Om\, {}^tx +
2x\,{}^tZ )}= e^{-\pi i\,\sigma (Z\Om^{-1}\,{}^tZ )}\cdot e^{\pi
i\,\sigma \{(x+Z\Om^{-1})\Om\,{}^t(x+Z\Om^{-1})\} }.$$ Since the
real Jacobi group $Sp(n,\BR)\ltimes H_\BR^{(m,n)}$ acts on
${\mathbb H}_{n,m}$ holomorphically, we may put
$$\Om=\,i\,A\,{}^t\!A,\quad Z=iV,\quad\  A\in\BR^{(n,n)},\quad
V=(v_{ij})\in\BR^{(m,n)}.$$
Then we obtain
\begin{eqnarray*}
& & \int_{\BR^{(m,n)}}  e^{\pi i\, \sigma ( x\Om\, {}^tx +
2x\,{}^tZ )} dx_{11}\cdots dx_{mn} \\
&=& e^{-\pi i\,\sigma (Z\Omega^{-1}\,{}^tZ)} \int_{\BR^{(m,n)}}
e^{\pi i\,\sigma [\{
x+iV(iA\,{}^t\!A)^{-1}\}(iA\,{}^t\!A)\,{}^t\!\{
x+iV(iA\,{}^t\!A)^{-1}\} ]}\,dx_{11}\cdots dx_{mn}\\
&=&e^{-\pi i\,\sigma (Z\Omega^{-1}\,{}^tZ)} \int_{\BR^{(m,n)}}
e^{\pi i\,\sigma [\{ x+V(A\,{}^t\!A)^{-1}\}A\,{}^t\!A\,{}^t\!\{
x+V(A\,{}^t\!A)^{-1}\} ]}\,dx_{11}\cdots dx_{mn}\\
&=& e^{-\pi i\,\sigma (Z\Omega^{-1}\,{}^tZ)} \int_{\BR^{(m,n)}}
e^{-\pi \,\sigma\{ (uA)\,{}^t\!(uA)\} }\,du_{11}\cdots
du_{mn}\quad \big(\,{\rm Put}\ u=
x+V(A\,{}^t\!A)^{-1}=(u_{ij}) \,\big)\\
&=& e^{-\pi i\,\sigma (Z\Omega^{-1}\,{}^tZ)} \int_{\BR^{(m,n)}}
e^{-\pi \,\sigma (w\,{}^t\!w)} (\det A)^{-m}\,dw_{11}\cdots
dw_{mn}\quad \big(\,{\rm Put}\ w=uA=(w_{ij})\,\big)\\
&=& e^{-\pi i\,\sigma (Z\Omega^{-1}\,{}^tZ)} \, (\det A)^{-m}\cdot
\left( \prod_{i=1}^m \prod_{j=1}^g \int_\BR e^{-\pi\,
w_{ij}^2}\,dw_{ij}\right)\\
&=& e^{-\pi i\,\sigma (Z\Omega^{-1}\,{}^tZ)} \, (\det A)^{-m}\quad
\big(\,{\rm because}\ \int_\BR e^{-\pi\,
w_{ij}^2}\,dw_{ij}=1\quad {\rm for\ all}\ i,j\,\big)\\
&=& e^{-\pi i\,\sigma (Z\Omega^{-1}\,{}^tZ)} \, \left( \det \big(
A\, {}^t\!A \big)\right)^{-{\frac m2}}\\
&=& e^{-\pi i\,\sigma (Z\Omega^{-1}\,{}^tZ)} \, \left( \det \left(
{ {\Omega}\over i } \right)\right)^{-{\frac m2}}.
\end{eqnarray*}

\noindent This completes the proof of Lemma 5.1. \hfill $\square$

\vskip 0.2cm According to Formula (4.8), for $x\in\rmn,$ we obtain

\begin{eqnarray*}
& & \left( \om_\CM \big(\wg \big) \mfoz\right) (x)\\
&=& \left( {\frac 1i}\right)^{{mn}\over 2}\big( \det
\CM\big)^{\frac n2}\,\int_{\rmn} \mfoz (y)\,e^{-2\pi\,i\,\s\,(\CM
\,y\,{}^tx)}dy\\
&=&\left( {\frac 1i}\right)^{{mn}\over 2}\big( \det
\CM\big)^{\frac n2}\,\int_{\rmn} e^{\pi\,i\,\sigma \{
\CM(y\,\Om\,{}^ty+2\, y\,{}^tZ)\}
}\,e^{-2\pi\,i\,\sigma(\CM \,y\,{}^tx)}dy\\
&=&\left( {\frac 1i}\right)^{{mn}\over 2}\big( \det
\CM\big)^{\frac n2}\,\int_{\rmn} e^{ \pi \,i\,\sigma\left\{
\CM\left( y\,\Om\,{}^ty\,+\,2\,y\,{}^t(Z-x)\right)\right\} }dy.
\end{eqnarray*}

\noindent If we substitute $u=\CM^{\frac 12}\,y,$ then $du=\left(
\det \CM\right)^{\frac n2}\,dy.$ Therefore according to Lemma 5.1,
we obtain

\begin{eqnarray*}
& & \left( \om_\CM \big(\wg \big) \mfoz\right) (x)\\
&=& \left( {\frac 1i}\right)^{{mn}\over 2}\big( \det
\CM\big)^{\frac n2}\,\int_{\rmn} e^{\pi\,i\,\sigma \left(
u\,\Om\,{}^tu\,+\,2\,\CM^{1/2}\,u\,{}^t(Z-x)\right)}\,\left(\det \CM\right)^{-{\frac n2}}du\\
&=&\left( {\frac 1i}\right)^{{mn}\over 2}\,\int_{\rmn}
e^{\pi\,i\,\sigma \left(
u\,\Om\,{}^tu\,+\,2\,u\,{}^t(\CM^{1/2}\,(Z-x))\right)}\,du\\
&=&\left( {\frac 1i}\right)^{{mn}\over 2}\, \left( \det
{\Omega\over i}\right)^{-{\frac m2}}\, e^{-\pi\,i\,\sigma\left\{
\CM^{1/2} (Z-x)\,\Om^{-1}\,{}^t(Z-x)\,\CM^{1/2}\right\}}\quad (\textrm{by\ Lemma\ 5.1})\\
&=&\,\left(\det \Om\right)^{-{\frac m2}}\,
e^{-\pi\,i\,\sigma\left( \CM\,(Z-x)\,\Om^{-1}\,{}^t(Z-x)\right)
}\\
&=&\,\left(\det \Om\right)^{-{\frac m2}}\,
e^{-\pi\,i\,\sigma\left(
\CM(Z\,\Om^{-1}\,{}^tZ\,+\,x\,\Om^{-1}\,{}^tx\,-\,2\,Z\,\Om^{-1}\,{}^tx)\right)
}.
\end{eqnarray*}

\noindent On the other hand, according to Formula (5.2), for $x\in
\rmn,$

\begin{eqnarray*}
& & J_\CM \big( \wg,(\Om,Z)\big)^{-1}\mfm_{\wg\cdot
(\Om,Z)}(x)\\
&=&\,e^{-\pi\,i\,\sigma(\CM\,Z\,\Om^{-1}\,{}^tZ)}\,\left(\det
\Om\right)^{-{\frac m2}}\, \mfm_{-\Om^{-1},Z\,\Om^{-1}}(x)\\
&=&\,\left(\det \Om\right)^{-{\frac
m2}}\,e^{-\pi\,i\,\sigma(\CM\,Z\,\Om^{-1}\,{}^tZ)}\,
e^{\pi\,i\,\sigma\left\{ \CM \left(
x\,(-\Om^{-1})\,{}^tx\,+\,2\,x\,{}^t(Z\,\Om^{-1}) \right)\right\}
}
\\
&=&\,\left(\det \Om\right)^{-{\frac m2}}\,
e^{-\pi\,i\,\sigma\left(
\CM(Z\,\Om^{-1}\,{}^tZ\,+\,x\,\Om^{-1}\,{}^tx\,-\,2\,Z\,\Om^{-1}\,{}^tx)\right)
}.
\end{eqnarray*}

\noindent Therefore we prove the covariance relation (5.3) in the
case $\wg=\sigma_n$. Since $J_\CM$ is an automorphic factor for
$G^J$ on $\BH_{n,m}$, we see that if the covariance relation (5.3)
holds for for two elements $\wg_1,\wg_2$ in $G^J$, then it holds
for $\wg_1 \wg_2.$ Finally we complete the proof. \hfill $\square$

\end{section}

\vskip 1cm

%
%

\begin{section}{{\bf Construction of Jacobi Forms }}
\setcounter{equation}{0}

\vskip 0.3cm Let $(\pi,V_\pi)$ be a unitary representation of
$G^J$ on the representation space $V_\pi$. We assume that
$(\pi,V_\pi)$ satisfies the following conditions (A) and (B):
\vskip 0.1cm \noindent {\bf (A)} There exists a vector valued map
\begin{equation*}
{\mathscr F}:\BH_{n,m}\lrt V_\pi,\quad\ (\Om,Z)\mapsto {\mathscr
F}_{\Om,Z}:={\mathscr F}(\Om,Z)
\end{equation*}

\noindent satisfying the following covariance relation
\begin{equation}
\pi\big( {\widetilde \g}\big) {\mathscr F}_{\Om,Z}=\psi\big(
{\widetilde \g}\big)\,J\big( {\widetilde
\g},(\Om,Z)\big)^{-1}\,{\mathscr F}_{ {\widetilde \g}\cdot (\Om,Z)
}\quad \textrm{for all}\ {\widetilde \g}\in G^J,\
(\Om,Z)\in\BH_{n,m},
\end{equation}

\noindent where $\psi$ is a character of $G^J$ and $J:G^J\times
\BH_{n,m}\lrt GL(1,\BC)$ is a certain automorphic factor for $G^J$
on $\BH_{n,m}.$

\vskip 0.1cm \noindent {\bf (B)} Let ${\widetilde \G}$ be an
arithmetic subgroup of $\G^J$. There exists a linear functional
$\theta:V_\pi\lrt \BC$ which is semi-invariant under the action of
${\widetilde \G}$, in other words, for all ${\widetilde\g}\in
{\widetilde\G}$ and $(\Om,Z)\in\BH_{n,m},$

\begin{equation}
\langle\, \pi^*\big( {\widetilde\g}\big)\theta,\,{\mathscr
F}_{\Om,Z}\,\rangle=\langle\, \theta,\pi\big(
{\widetilde\g}\big)^{-1}{\mathscr
F}_{\Om,Z}\,\rangle=\chi\big({\widetilde\g}\big)\,\langle\,
\theta,\,{\mathscr F}_{\Om,Z}\,\rangle ,
\end{equation}

\noindent where $\pi^*$ is the contragredient of $\pi$ and
$\chi:{\widetilde\G}\lrt T$ is a unitary character of
${\widetilde\G}$.

\vskip 0.2cm Under the assumptions (A) and (B) on a unitary
representation $(\pi,V_\pi)$, we define the function $\Theta$ on
$\BH_{n,m}$ by

\begin{equation}
\Theta(\Om,Z):=\,\langle\,\theta,{\mathscr
F}_{\Om,Z}\,\rangle=\theta\big({\mathscr F}_{\Om,Z}\big),\quad\
(\Om,Z)\in\BH_{n,m}.
\end{equation}

We now shall see that $\Theta$ is an automorphic form on
$\BH_{n,m}$ with respect to ${\widetilde\G}$ for the automorphic
factor $J$.

\begin{lemma} Let $(\pi,V_\pi)$ be a unitary representation of
$G^J$ satisfying the above assumptions (A) and (B). Then the
function $\Theta$ on $\BH_{n,m}$ defined by (6.3) satisfies the
following modular transformation behavior

\begin{equation}
\Theta\big({\widetilde \g}\cdot(\Om,Z)\big)=\,\psi\big(
{\widetilde
\g}\big)^{-1}\,\chi\big({\widetilde\g}\big)^{-1}\,J\big(
{\widetilde \g},(\Om,Z)\big)\,\Theta(\Om,Z)
\end{equation}

\noindent for all ${\widetilde\g}\in {\widetilde\G}$ and
$(\Om,Z)\in\BH_{n,m}.$
\end{lemma}

\noindent {\it Proof.} For any ${\widetilde\g}\in {\widetilde\G}$
and $(\Om,Z)\in\BH_{n,m},$ according to the assumptions (6.1) and
(6.2), we obtain

\begin{eqnarray*}
& & \Theta\big({\widetilde
\g}\cdot(\Om,Z)\big)=\big\langle\,\theta,
{\mathscr F}_{ {\widetilde \g}\cdot (\Om,Z)} \big\rangle \hskip 5cm\\
&=&\big\langle\,\theta,\psi\big( {\widetilde \g}\big)^{-1}J\big(
{\widetilde \g},(\Om,Z)\big)\,\pi\big({\widetilde\g}\big){\mathscr
F}_{ \Om,Z}\,\big\rangle\\
&=&\psi\big( {\widetilde \g}\big)^{-1} J\big( {\widetilde
\g},(\Om,Z)\big)\,\big\langle\,\theta,\pi\big({\widetilde\g}\big){\mathscr
F}_{ \Om,Z}\,\big\rangle\\
&=&\,\psi\big( {\widetilde
\g}\big)^{-1}\,\chi\big({\widetilde\g}\big)^{-1}\,J\big(
{\widetilde \g},(\Om,Z)\big)\,\big\langle\,\theta,{\mathscr
F}_{ \Om,Z}\,\big\rangle\\
&=&\,\psi\big( {\widetilde
\g}\big)^{-1}\,\chi\big({\widetilde\g}\big)^{-1}\,J\big(
{\widetilde \g},(\Om,Z)\big)\,\Theta(\Om,Z).
\end{eqnarray*}
\hfill$\square$

\newcommand\zmn{\BZ^{(m,n)} }
\newcommand\wgam{\widetilde\gamma}

\vskip 0.2cm Now for a positive definite integral symmetric matrix
$\CM$ of degree $m$, we define the holomorphic function
$\Theta_\CM:\BH_{n,m}\lrt\BC$ by

\begin{equation}
\Theta_\CM (\Om,Z):=\sum_{\xi\in \BZ^{(m,n)}}
e^{\pi\,i\,\sigma\left( \CM(
\xi\,\Om\,{}^t\xi\,+\,2\,\xi\,{}^tZ)\right) },\quad (\Om,Z)\in
\BH_{n,m}.
\end{equation}

\begin{theorem} Let $\CM$ be a symmetric positive definite, unimodular even integral matrix of degree $m$.
Then for any $\wgam=(\g,(\l,\mu;\kappa))\in\Gamma^J$ with $\g \in
\Gamma_n$ and $(\l,\mu;\kappa)\in H_{\BZ}^{(n,m)},$ the function
$\Theta_\CM$ satisfies the functional equation
\begin{equation}
\Theta_\CM\big( \wgam\cdot (\Om,Z)\big)=\rho_\CM (\wgam)\,
J_\CM\big( \wgam,(\Om,Z)\big) \Theta_\CM(\Om,Z),\quad (\Om,Z)\in
\BH_{n,m},
\end{equation}
\end{theorem}

\noindent where $\rho_\CM(\wgam)$ is a uniquely determined character of $\G^J$ with
$|\rho_\CM(\wgam)|^8=1$ and $J_\CM:G^J\times \BH_{n,m}\lrt
\BC^{\times}$ is the automorphic factor for $G^J$ on $\BH_{n,m}$
defined by the formula (5.2).

\noindent {\it Proof.}  For an element
$\wgam=(\g,(\lambda,\mu;\kappa))\in \G^J$ with $\g=\begin{pmatrix}
A & B
\\ C & D \end{pmatrix}\in \G_{n}$ and $(\l,\mu;\kappa)\in H_\BZ^{(n,m)},$ we put
$(\Om_*,Z_*)=\wgam\cdot (\Om,Z)$ for $(\Om,Z)\in\BH_{n,m}.$ Then
we have
\begin{eqnarray*}
&\Om_*=\g\cdot \Om=(A\Om+B)(C\Om+D)^{-1},\\
&Z_*=(Z+\l\, \Om+\mu)(C\Om+D)^{-1}.
\end{eqnarray*}

\noindent We define the linear functional $\vartheta$ on $L^2\big(
\rmn\big)$ by

\begin{equation*}
\vartheta (f)=\langle \vartheta,f \rangle:=\sum_{\xi\in
\zmn}f(\xi),\quad\ f\in L^2\big( \rmn\big).
\end{equation*}

\noindent We note that $\Theta_\CM(\Om,Z)=\vartheta\big(
\mfoz\big).$ Since $\mfm$ is a covariant map for the
Schr{\"o}dinger-Weil representation $\omega_\CM$ by Theorem 5.1,
according to Lemma 6.1, it suffices to prove that $\vartheta$ is
semi-invariant for $\om_\CM$ under the action of $\G^J$, in other
words, $\vartheta$ satisfies the following semi-invariance
relation
\begin{equation}
\Big\langle\, \vartheta,\om_\CM\big( \wgam\big)\mfoz \,\Big\rangle
=\,\rho_\CM \big( \wgam\big)^{-1} \,\Big\langle\, \vartheta,\mfoz
\,\Big\rangle
\end{equation}

\noindent for all $\wgam\in \G^J$ and $(\Om,Z)\in \BH_{n,m}.$

\vskip 0.12cm We see that the following elements
$h(\lambda,\mu;\kappa),\ t(b),\,g(\alpha)$ and $\s_n$ of $\G^J$
defined by
\begin{eqnarray*}
&& h(\la,\mu;\kappa)=\big( I_{2n},(\la,\mu;\kappa)\big)\
\textrm{with}\  \la,\mu\in
\BZ^{(m,n)}\ \textrm{and}\ \kappa\in\BZ^{(m,m)} ,\\
&&t(b)=\big( t_0(b),(0,0;0)\big)\ \textrm{with any}\
                   b=\,{}^tb\in \BZ^{(n,n)},\\
&& g(\alpha)=\big( g_0(\alpha),(0,0;0)\big)\
\textrm{with any}\ \alpha\in GL(n,\BZ),\\
 &&\s_n=\big( s_{n,0},(0,0;0)\big)
\end{eqnarray*}
generate the Jacobi modular group $\G^J.$ Therefore it suffices to
prove the semi-invariance relation (6.7) for the above generators
of $\G^J.$

\vskip 0.5cm\noindent {\bf Case I.} $\wgam=h(\l,\mu;\kappa)$ with
$\l,\mu\in\zmn,\ \kappa\in\BZ^{(m,m)}.$

\vskip 0.1cm In this case, we have
$$\Om_*=\Om,\quad Z_*=Z+\la\,\Om+\mu$$
and
$$J_\CM\big(\wgam,(\Om,Z)\big)=\,e^{-\pi i\,\s\{ \CM(\l\,\Om\,{}^t\l
+2\,\l\,{}^tZ+\kappa+\mu\,{}^t\!\l)\} }.$$

According to the covariance relation (5.3),
\begin{eqnarray*}
& &\big\langle \,\vartheta, \om_\CM\big( \wgam\big)
\mfoz\,\big\rangle\\
&=&\,\big\langle \,\vartheta,
J_\CM\big( \wgam,(\Om,Z)\big)^{-1}\mfm_{\wgam\cdot(\Om,Z)}\,\big\rangle\\
&=&\,J_\CM\big( \wgam,(\Om,Z)\big)^{-1}\,\big\langle \,\vartheta,
\mfm_{\Om,Z+\la\,\Om+\mu}\,\big\rangle\\
&=&\,J_\CM\big(
\wgam,(\Om,Z)\big)^{-1}\sum_{A\in\zmn}e^{\pi\,i\,\s\left\{ \CM
\left( A\Om\,{}^t\!A\,+\,2\,A\,{}^t(Z+\l\,\Om+\mu)\right) \right\}
}\\
&=&\,J_\CM\big( \wgam,(\Om,Z)\big)^{-1}\cdot
e^{-\pi\,i\,\sigma\left( \CM
(\l\,\Om\,{}^t\!\l\,+\,2\,\l\,{}^tZ)\right)}\\
& &\times \sum_{A\in\zmn} e^{2\,\pi\,i\,\sigma(\CM A\,{}^t\mu)}
e^{\pi\,i\,\s\left\{ \CM \left(
(A+\l)\,\Om\,{}^t\!(A+\l)\,+\,2\,(A+\l)\,{}^tZ\right) \right\}
}\\
&=& e^{\pi\,i\,\sigma\left( \CM(\kappa\,+\,\mu\,{}^t\!\l)\right)}
\,\big\langle \,\vartheta, \mfm_{\Om,Z}\,\big\rangle.
\end{eqnarray*}

\noindent Here we used the fact that $\s(\CM A\,{}^t\mu)$ is an
integer. We put $\rho_\CM\big( \wgam\big)=\rho_\CM\big(
h(\l,\mu;\kappa)\big)=e^{-\pi\,i\,\sigma\left( \CM
(\kappa\,+\,\mu\,{}^t\!\l)\right)}.$ Therefore $\vartheta$
satisfies the semi-invariance relation (6.7) in the case
$\wgam=h(\l,\mu;\kappa)$ with $\l,\mu\in\zmn,\
\kappa\in\BZ^{(m,m)}.$

\vskip 0.5cm\noindent {\bf Case II.} $\wgam=t(b)$ with
$b=\,{}^tb\in\BZ^{(n,n)}.$

\vskip 0.1cm In this case, we have
$$\Om_*=\Om+b,\quad Z_*=Z\quad \textrm{and}\quad J_\CM\big(\wgam,(\Om,Z)\big)=1.$$

\noindent According to the covariance relation (5.3), we obtain
\begin{eqnarray*}
& &\big\langle \,\vartheta, \om_\CM\big( \wgam\big)
\mfoz\,\big\rangle\\
&=&\,\big\langle \,\vartheta,
J_\CM\big( \wgam,(\Om,Z)\big)^{-1}\mfm_{\wgam\cdot(\Om,Z)}\,\big\rangle\\
&=&\,\big\langle \,\vartheta,
\mfm_{\Om+b,Z}\,\big\rangle\\
&=&\, \sum_{A\in\zmn} e^{\pi\,i\,\s\left\{ \CM \left(
A\,(\Om+b)\,{}^t\!A\,+\,2\,A\,{}^tZ\right) \right\}
}\\
&=&\, \sum_{A\in\zmn} e^{\pi\,i\,\s\left( \CM \left(
A\,\Om\,{}^t\!A\,+\,2\,A\,{}^tZ\right) \right) }\cdot
e^{\pi\,i\,\sigma(\CM A\,b\,{}^t\!A)}\\
&=&\, \sum_{A\in\zmn} e^{\pi\,i\,\sigma\left( \CM (
A\,\Om\,{}^t\!A\,+\,2\,A\,{}^tZ ) \right) }\\
&=&\,\big\langle \,\vartheta, \mfm_{\Om,Z}\,\big\rangle.
\end{eqnarray*}

\noindent Here we used the fact that $\s(\CM A\,b\,{}^t\!A)$ is an
even integer. We put $\rho_\CM\big( \wgam\big)=\rho_\CM\big(
t(b)\big)=1.$ Therefore $\vartheta$ satisfies the semi-invariance
relation (6.7) in the case $\wgam=t(b)$ with
$b=\,{}^tb\in\BZ^{(n,n)}.$

\vskip 0.5cm\noindent {\bf Case III.} $\wgam=g(\alpha)$ with
$\alpha\in GL(n,\BZ).$

\vskip 0.1cm In this case, we have
$$\Om_*=\,{}^t\alpha\,\Om\,\alpha,\quad Z_*=Z\alpha$$
and
$$J_\CM\big(\wgam,(\Om,Z)\big)=(\det\alpha)^{-{\frac m2}}.$$

\noindent According to the covariance relation (5.3), we obtain
\begin{eqnarray*}
& &\big\langle \,\vartheta, \om_\CM\big( \wgam\big)
\mfoz\,\big\rangle\\
&=&\,\big\langle \,\vartheta,
J_\CM\big( \wgam,(\Om,Z)\big)^{-1}\mfm_{\wgam\cdot(\Om,Z)}\,\big\rangle\\
&=&\,\left( \det \alpha \right)^{\frac m2}\,\big\langle
\,\vartheta,
\mfm_{{}^t\alpha\,\Om\,\alpha,Z\,\alpha}\,\big\rangle\\
&=&\,\left( \det \alpha \right)^{\frac m2} \sum_{A\in\BZ^{(m,n)} }
\mfm_{{}^t\alpha\,\Om\,\alpha,Z\,\alpha}(A)\\
&=&\,\left( \det \alpha \right)^{\frac m2}\sum_{A\in\BZ^{(m,n)} }
e^{\pi i\,\sigma\{
\CM(A\,{}^t\alpha\,\Om\,{}^t(A\,{}^t\alpha)+2\,A\,{}^t\alpha\,{}^tZ)\}} \\
&=&\,\left( \det \alpha \right)^{\frac m2}\,\big\langle
\,\vartheta, \mfm_{\Om,Z}\,\big\rangle.
\end{eqnarray*}

\noindent Here we put $\rho_\CM\big( \wgam\big)=\rho_\CM\big(
g(\alpha)\big)=(\det\alpha)^{-{\frac m2}}.$ Therefore $\vartheta$
satisfies the semi-invariance relation (6.7) in the case
$\wgam=g(\alpha)$ with $\alpha\in GL(n,\BZ).$

\vskip 0.52cm\noindent {\bf Case IV.} $\wgam=\left(
\begin{pmatrix} 0 & -I_n \\ I_n & \ 0
\end{pmatrix},(0,0;0)\right).$

\vskip 0.1cm In this case, we have
$$\Om_*=-\Om^{-1},\quad Z_*=Z\,\Om^{-1}$$
and
$$J_\CM\big(\wgam,(\Om,Z)\big)=e^{\pi \,i\,\sigma (\CM Z\Om^{-1}\,{}^tZ)}\,\big( \det\Om\big)^{\frac m2}.$$

\noindent In the process of the proof of Theorem 5.1, using Lemma
5.1, we already showed that

\begin{equation}
\int_{\rmn} e^{\pi\,i\,\sigma( \CM(y\,\Om\,{}^ty\,+\,2\,y\,{}^tZ)
) } dy =\,\big(\det \CM\big)^{-{\frac n2}}\left( \det {\Om \over
i}\right)^{-{\frac
m2}}\,e^{-\pi\,i\,\sigma(\CM\,Z\,\Om^{-1}\,{}^tZ)}.
\end{equation}

\noindent By (6.8), we see that the Fourier transform of $\mfoz$
is given by

\begin{equation}
\widehat {\mfoz} (x)=\,\big(\det \CM\big)^{-{\frac n2}}\left( \det
{\Om \over i}\right)^{-{\frac
m2}}\,e^{-\pi\,i\,\sigma(\CM\,(Z-x)\,\Om^{-1}\,{}^t(Z-x))}.
\end{equation}

\noindent According to the covariance relation (5.3), Formula
(6.9) and Poisson summation formula, we obtain

\begin{eqnarray*}
& &\big\langle \,\vartheta, \om_\CM\big( \wgam\big)
\mfoz\,\big\rangle\\
&=&\,\big\langle \,\vartheta,
J_\CM\big( \wgam,(\Om,Z)\big)^{-1}\mfm_{\wgam\cdot(\Om,Z)}\,\big\rangle\\
&=&\,J_\CM\big( \wgam,(\Om,Z)\big)^{-1} \big\langle \,\vartheta,
\mfm_{-\Om^{-1},Z\Om^{-1}}\,\big\rangle\\
&=&\,(\det \Om)^{-{\frac
m2}}\,e^{-\pi\,i\,\sigma(\CM\,Z\,\Om^{-1}\,{}^tZ) }
\sum_{A\in\BZ^{(m,n)}} e^{-\pi\,i\,\sigma\left( \CM(
A\,\Om^{-1}\,{}^t\!A\,-\,2\,A\,\Om^{-1}\,{}^tZ) \right)}\\
&=&\,(\det \Om)^{-{\frac m2}}\, \sum_{A\in\BZ^{(m,n)}}
e^{-\pi\,i\,\sigma\left( \CM(Z\,\Om^{-1}\,{}^tZ\,+\,
A\,\Om^{-1}\,{}^t\!A\,-\,2\,A\,\Om^{-1}\,{}^tZ) \right)}\\
&=&\,(\det \Om)^{-{\frac m2}}\, \sum_{A\in\BZ^{(m,n)}}
e^{-\pi\,i\,\sigma\left( \CM(Z-A)\,\Om^{-1}\,{}^t(Z-A)\right)}\\
&=&\,(\det \Om)^{-{\frac m2}}\big(\det \CM\big)^{{\frac n2}}\left(
\det {\Om \over i}\right)^{{\frac m2}}\,\sum_{A\in\BZ^{(m,n)}}
\widehat {\mfoz} (A) \quad (\,\textrm{by\ Formula}\ (6.9))\\
&=&\,\big(\det \CM\big)^{{\frac n2}}\left( \det {{I_n} \over
i}\right)^{{\frac m2}}\,\sum_{A\in\BZ^{(m,n)}}
{\mfoz} (A) \quad (\,\textrm{by\ Poisson summation formula}) \\
&=&\,\big(\det \CM\big)^{{\frac n2}}\,(-i)^{{mn}\over
2}\,\big\langle \,\vartheta, \mfm_{\Om,Z}\,\big\rangle \\
&=&\,(-i)^{{mn}\over 2}\,\big\langle \,\vartheta,
\mfm_{\Om,Z}\,\big\rangle.
\end{eqnarray*}

\noindent Here we used the fact that $\det \CM=1$ because $\CM$ is
unimodular. We put $\rho_\CM\big( \wgam\big)=\rho_\CM (\sigma_n)
=(-i)^{-{{mn}\over 2}}.$ Therefore $\vartheta$ satisfies the
semi-invariance relation (6.7) in the case $\wgam=\sigma_n.$ The
proof of Case IV is completed. Since $J_\CM$ is an automorphic
factor for $G^J$ on $\BH_{n,m}$, we see that if the formula (6.6)
holds for two elements $\wgam_1,\wgam_2$ in $\G^J$, then it holds
for $\wgam_1 \wgam_2.$ Finally we complete the proof of Theorem
6.1. \hfill $\square$

\vskip 0.2cm
\begin{remark}
For a symmetric positive
definite integral matrix $\CM$ that is not unimodular even
integral, we obtain a similar transformation formula like (6.6).
If $m$ is odd, $\Theta_\CM(\Om,Z)$ is a Jacobi form of a
half-integral weight ${\frac m2}$ and index ${\CM}\over 2$ with
respect to a suitable arithmetic subgroup $\G^J_{\Theta,\CM}$ of
$\G^J$ and a character $\rho_\CM$ of $\G^J_{\Theta,\CM}$.
\end{remark}
 
\vskip 0.2cm For instance, we obtain the following\,:

\begin{theorem}
Let $\CM$ be a symmetric positive definite integral matrix of degree $m$
such that $\det\, (\CM)=\,1$.
Let $\Gamma_{1,2}$ be an arithmetic subgroup of $\G_n$ generated by all the following elements
\begin{equation*}
t(b)=\, \begin{pmatrix} I_n & b \\ 0 & I_n \end{pmatrix},\quad
g(\alpha)=\,\begin{pmatrix} {}^t\alpha & 0 \\ 0 & \alpha^{-1} \end{pmatrix}, \quad
\sigma_n=\,\begin{pmatrix} 0 & -I_n \\ I_n & \ \ 0 \end{pmatrix},
\end{equation*}
where $b=\,{}^tb\in\BZ^{(n,n)}$ with {\it even\ diagonal} and $\alpha\in \BZ^{(n,n)}.$
We put
$$\G_{1,2}^J:=\,\G_{1,2}\ltimes H_\BZ^{(n,m)}.$$
Then $\Theta_\CM$ satisfies the transformation formula (6.6) for all 
${\widetilde \g}\in \G_{1,2}^J.$ Therefore $\Theta_\CM$ is a Jacobi form of weight ${\frac m2}$
with level $\G_{1,2}$ and index ${\CM}\over 2$ for the uniquely determined character $\rho_\CM$
of $\G_{1,2}^J$.
\end{theorem}

\vskip 0.3cm \noindent
{\it Proof.} The proof is essentially the same as the proof of Theorem 6.1. We leave the detail
to the reader.
\hfill $\square$

\end{section}


\vskip 2cm
\bibliography{central}

\end{document}